\documentclass{amsart} 

\usepackage{epsf, graphicx} 
\usepackage{amssymb,latexsym, amsmath, amscd, array} 
\def\hepsffile{\leavevmode\epsffile} 

\swapnumbers 
\numberwithin{equation}{section} 

\theoremstyle{plain} 
\newtheorem{thm}{Theorem}[section] 
\newtheorem{cor}[thm]{Corollary}

\newtheorem{prop}[thm]{Proposition}

\newcommand\propref{Proposition~\ref}

\theoremstyle{definition} 
\newtheorem{defin}[thm]{Definition}

\newtheorem{rem}[thm]{Remark}

\def\id{\protect\operatorname{id}}

\def\pr{\protect\operatorname{proj}}

\def\nn{\protect\operatorname{\mathcal N}} 
\def\NN{\protect\operatorname{\mathcal N}} 
\def\PP{\protect\operatorname{\mathcal P}} 
\def\li{\protect\operatorname{lift}}

\def\ga{\alpha} 
 
\def\gg{\gamma}

\def\gf{\varphi} 
\def\bor{\Omega} 
\def\om{\omega}

\def\Z{{\mathbb Z}}

\def\N{{\mathbb N}} 
\def\oma{\mathbb O} 
\def\n0{{\mathbb N\cup\{0\}}}

\def\1{\hbox{\rm\rlap {1}\hskip.03in{\rom I}}} 
\def\Bbbone{{\rm1\mathchoice{\kern-0.25em}{\kern-0.25em} 
{\kern-0.2em}{\kern-0.2em}I}} 

\def\empt{\varnothing}

\def\pp{\medskip{\parindent 0pt \it Proof.\ }}

\def\m{\medskip} 
\def\ov{\overline} 
 
\def\fm{\mathfrak m} 
\def\fmm{\mathfrak M}

\begin{document} 
\date{June 8, 2003} 

\leftline{ } 
\centerline{ } 
\title[Algebraic structures on generalized strings] 
{Algebraic structures on generalized strings} 
\author[V.~Chernov (Tchernov) and Yu.~Rudyak]{Vladimir V. 
Chernov 
(Tchernov) and Yuli B. Rudyak} 
\address{V. Chernov, Department of Mathematics, 
6188 Bradley Hall, Dartmouth College, Hanover NH 03755, 
USA} 
\email{chernov@math.ethz.ch} 
\address{Yu. Rudyak, Department of Mathematics, University 
of Florida, 358 
Little Hall, Gainesville, FL 32611-8105, USA} 
\email{rudyak@math.ufl.edu} 

\Large 

\begin{abstract} A garland based on a manifold $P$ is a finite set of 
manifolds homeomorphic to $P$ with some of them glued together at 
marked points. Fix a manifold $M$ and consider a space $\NN$ of all 
smooth mappings of garlands based on $P$ into $M$. We construct 
operations $\bullet$ and $[-,-]$ on the bordism groups $\bor_*(\NN)$ 
that give $\bor_*(\NN)$ the natural graded commutative assosiative  and 
graded Lie algebra structures. We also construct two auto-homomorphisms 
$\pr$ and $\li$ of $\bor_*(\NN)$ such that 
$\pr(\li \alpha_1\bullet \li \alpha_2)= [\alpha_1, \alpha_2]$ for all 
$\alpha_1, \alpha_2 \in \bor_*(\NN )$.
If $P$ is a boundary, then $\pr \circ \li=0$ and thus $\Delta^2=0$ 
for $\Delta=\li \circ \pr$. We show that under certain conditions the 
operations $\Delta$ and $\bullet$ give rise to Batalin-Vilkoviski and 
Gerstenhaber algebra structures on $\bor_*(\NN)$.

In a particular case when $P=S^1$, the algebra $\bor_*(\NN)$ is 
related to
the string-homology algebra constructed by Chas and 
Sullivan~\cite{ChasSullivan}.
\end{abstract}
\maketitle 

\section{Preliminaries} 
We work in the smooth category.

Given a topological space $X$, we define a {\it mark} $\fm$ on $X$ to be 
a 
finite family $q_{\gamma}, \gamma \in \Gamma,$ of points of $X$, the 
case $\Gamma=\empt$ is allowed. It can also happen that $q_{\gg}= 
q_{\gamma'}$ for $\gg\ne \gg'$. The number of elements in a mark $\fm$ 
is denoted by $|\fm|$ or $|\Gamma|$. 

\m
Given a marked space $X=(X,\fm)$, a {\it marked map} $f: X \to Y$ is a 
map such that $f(q_{\gg})=f(q_{\gg'})$ for all $\gg,\gg'\in \Gamma$. 
Frequently, if we need to indicate the mark $\fm$ on $X$, we will write 
$f: (X,\fm) \to Y$ and say that $f$ is an $\fm$-marked map. 

\m
For future goals, we assign to each mark $\fm$ its {\it grading} $g(\fm) 
\in \N$. 

\m A multimark $\fmm$ on a space $X$ is just a finite family 
$\{\fm_1,\ldots , \fm_k\}$ of marks on $X$. It can happen that the 
different marks can have common elements, or even $\fm_i=\fm_j$ for 
$i\ne j$. A multimarked map $f: (X,\fmm) \to Y$ is just a map which is a 
marked map for each mark $\fm \in \fmm$. 

\m Fix a closed connected oriented manifold $P$ of dimension $n$ and a 
connected oriented manifold $M$ of dimension $m$, and let $\PP$ be the 
set of all smooth maps $P \to M$ in $C^{\infty}$-topology. Fix $k\in 
\N$, $l\in \n0$ and $g_1, \ldots, g_l\in \N$. 
Consider a multimarked manifold $(P_1\sqcup \cdots \sqcup P_k; \fm_1, 
\ldots, \fm_l)$ with $g(\fm_i)=g_i$ for $i=1, \ldots, l$. Then the 
set $\nn(k;l,g_1, \ldots , g_l)$ of all multimarked maps $(P_1\sqcup 
\cdots \sqcup P_k; \fm_1, \ldots, \fm_l) \to M$ factorized by the action 
of the permutation group $S_k$ can be regarded as a subspace of the 
space $\bigl(P^{\sum |\fm_i|}\times \PP^k \bigr) /S_k$, and we equip 
$\nn(k;l,g_1, \ldots , g_l)$ with the subspace topology. 

\m Now, we set 
$$ 
\nn=M\sqcup 
\coprod_{k=1}^{\infty}\coprod_{l=0}^{\infty}\coprod_{g_1, \ldots, 
g_l}\nn(k;l, g_1, \ldots , g_l). 
$$ 
Here the first summand $M$ can be treated as $\nn(0;0)$. 

\begin{rem}\label{garland} 
Given a multimarked manifold $(P\sqcup \cdots \sqcup P, \fmm)/S_k\in 
\PP_k/S_k$, consider the following equivalence relation $\simeq$ on it: 
Two points $a,b\in P\sqcup \cdots \sqcup P$ are equivalent if there 
exists a mark $\fm \in \fmm$ such that $a,b\in \fm$. We define a {\it 
garland $($based on the manifold $P)$} of $(P\sqcup \cdots \sqcup P, 
\fmm)$ to be the quotient space $\Bigl(\bigl(P\sqcup \cdots \sqcup 
P)/\simeq\bigr)/S_k\Bigr)$, see~Figure~\ref{garland.fig}. 

\begin{figure}[htbp] 
\begin{center} 
\epsfxsize\hsize\advance\epsfxsize -0.5cm 
\epsfxsize 10 cm 
\hepsffile{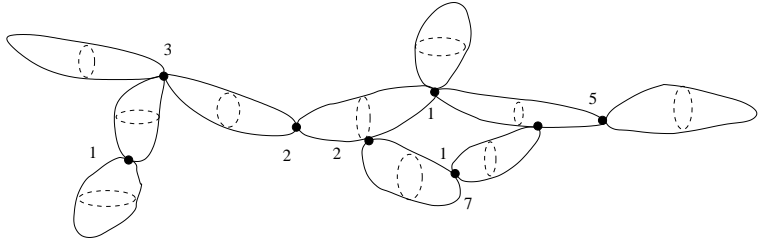} 
\end{center} 
\caption{A garland with $\N$-graded marks}\label{garland.fig} 
\end{figure}

Clearly, every space $\nn_k, k>0,$ consists of maps of garlands. In 
fact, garlands are precisely the objects we will work with. However, 
since garlands are not manifolds and therefore the problems with 
transversality can appear, we prefer to interpret garlands as 
multimarked manifolds and deal with the last ones. 
\end{rem} 

\m Throughout the paper we work with oriented bordism theory 
$\bor_*(-)$. All the necessary information can be found in books 
\cite{Rudyak, Stong, Switzer}.

\section{An operation $\bullet$ on $\nn$}

\m In this section we construct a commutative and associative operation 
\begin{equation}\label{bullet} 
\bullet:\bor_i(\nn)\otimes \bor_j(\nn) \to \bor_{i+j-m}(\nn) 
\end{equation} 

\m 
Let $\ov\ga_1:F_1^i\to \nn$ and 
$\ov\ga_2:F_2^j\to \nn$ be representatives of $\ga_1\in \bor_i(\NN)$ and 
$\ga_2\in \bor_j(\NN)$, respectively. Without loss of generality we 
assume that $F_1^i, F_2^j$ are connected. So, we have adjoint maps 
$\om_1: F_1\times N_1 \to M$ and $\om_2: F_2\times N_2 \to M$. If, say, 
$\ov\ga_2(F_2^j)\subset M=N(0;0)\subset \nn$ then $\om_2$ is just 
$\ov\ga_2:F_2 \to M$. 

\m Let $\fm_1, \fm_2$ be two marks of grading 1 on $N_1$ and $N_2$, 
respectively. Choose $q_1\in \fm_1$ and $q_2\in \fm_2$ and 
consider the pull-back diagram 
\begin{equation}\label{v-bordism} 
\CD 
V_{q_1, q_2}@>j_1>> F^i_1\times q_1\\ 
@VVj_2V @VV\omega _1 V\\ 
F^j_2\times q_2 @>\omega _2 >> M 
\endCD 
\end{equation} 
(if $\ov\ga_2(F_2^j)\subset M=N(0;0)\subset \nn$ then $F^j_2\times 
q_2$ is just $F^j_2$). Using standard transversality arguments, we can 
and shall assume that 
$V=V_{q_1, q_2}$ is a smooth $(i+j-m)$-dimensional oriented manifold. 

Together the maps 
\begin{equation*} 
\CD 
V\times N_1 @>j_1\times \id >> F_1\times N_1 @>\om_1>> M 
\endCD 
\end{equation*} 
and 
\begin{equation*} 
\CD 
V\times N_2 @>j_2\times \id >> F_2\times N_2 @>\om_1>> M 
\endCD 
\end{equation*} 
yield the map $\gf: V \times (N_1\sqcup N_2) \to M$. 
We equip $N_1\sqcup N_2$ with the mark $\fm_1\cup \fm_2$ and assign to 
this mark grading 1. 
All the other marks on $N_1$ and $N_2$ remain as they were and do not 
change their gradings. 
Clearly, for every $v\in V$, the map 
\begin{equation*} 
\CD 
N_1\sqcup N_2 = \{v\}\times (N_1\sqcup N_2) @>\gf >> M 
\endCD 
\end{equation*} 
is a marked map. So, the image of the adjoint to $\gf$ map 
$\psi=\psi_{q_1,q_2}$  belongs to $\nn$, i.e. we have a map 
$\psi_{q_1,q_2}: V \to \nn$. 

It is easy to see that the bordism class of 
$$ 
\bigl [ V=V_{\fm _1,\fm _2}, 
\psi=\psi_{\fm _1,\fm _2}\bigr ] 
$$ 
depends only on $\fm_1, \fm_2$ (and not on the choices of $q_1\in \fm_1, 
q_2\in \fm_2$) and on the bordism classes 
$\ga_1, \ga_2$. Now, we define $\ga_1\bullet \ga_2=\bigl [ V=V_{\fm 
_1,\fm _2}, 
\psi=\psi_{\fm _1,\fm_2}\bigr ]$.

\m Put 
%
%
\begin{equation}\label{bulletdefin} 
\ga_1\bullet \ga_2=\sum_{\fm_1, \fm_2}[V_{\fm _1, 
\fm _2}, \psi_{\fm _1, \fm _2}]\in \bor _{i+j-m}(\nn), 
\end{equation} 

where $\fm_1$ and $\fm_2$ run over all marks of grading one on $N_1$ and 
$N_2$, respectively. 

\m Following \cite{ChasSullivan}, we set 
$\oma_i(\nn)=\bor_{i+m}(\nn)$. Now the pairing $\bullet$ gets the form 

\begin{equation}\label{omabullet} 
\bullet: \oma_i (\NN )\otimes \oma_j(\NN )\to \oma _{i+j}(\NN) 
\end{equation} 

We say that $|\ga|=i$ whenever $\ga\in \oma_i(\nn)$. The following 
theorem follows directly from the definition of the operation $\bullet$.

\begin{thm}\label{bulletasso} 
The operation 
$\bullet: \oma_i (\NN )\otimes \oma_j(\NN )\to \oma _{i+j}(\NN)$, 
$i,j\in 
\Z$, 
converts the abelian group $\oma_*(\NN)$ into an associative graded 
commutative ring, i.e. the operation $\bullet$ has the following 
properties: 
\begin{description} 
\item[1] $\ga_1 \bullet \ga _2 = (-1)^{|\ga_1||\ga_2|}\ga_2 \bullet 
\ga _1$, for 
all $\ga_1, \ga_2\in \oma_*(\NN)$; 
\item[2] $(\ga_1 \bullet \ga_2)\bullet \ga _3=\ga_1 \bullet 
(\ga_2 \bullet \ga _3)$, for all $\ga_1, \ga_2, \ga_3 
\in \oma_*(\NN)$; 
\item[3] $\ga_1\bullet (\ga_2 + \ga _3)=\ga_1\bullet \ga _2 
+\ga_1 \bullet \ga_3$ for all $\ga_1, \ga_2, \ga_3\in 
\oma_*(\NN)$. 
\end{description} 
\qed 
\end{thm} 

We leave it to the reader to check that the unit of the operation 
$\bullet$ is given by the bordism class of the inclusion 
$M=\nn(0;0)\to \NN$.

\section{The string bracket $[-,-]$ on $\nn$}\label{sectionbracket} 

Let all the notation be as in the previous section. Consider the 
pull-back diagram 
\begin{equation}\label{w-bordism} 
\CD 
W @>j_1>> F^i_1\times N_1\\ 
@VVj_2V @VV\ \omega _1 V\\ 
F^j_2\times N_2 @> \omega _2 >> M\\ 
\endCD 
\end{equation} 

Using standard transversality arguments, we can and shall assume that 
$W$ is a smooth $(i+j+2n-m)$-dimensional oriented manifold. Define the 
map 

\begin{equation*} 
\CD 
a_k: W @>j_k >> F_k\times N_k @>p_2 >> N_k,\, @. k=1,2. 
\endCD 
\end{equation*} 

For every $w\in W$, we equip the manifold $N_1\sqcup N_2$ with the mark 
$\fm_w=\{a_1(w), a_2(w)\}$ which, by definition, has grading 2. All 
the other marks on $N_1$ and $N_2$ remain as they were and do not change 
their gradings. Now, for every $w\in W$ we have a marked map $\psi_w: 
N_1\sqcup N_2 \to M$ where, say, for $n\in N_1$ we have 
$\psi_w(n)=\om_1(j_1(w), n)$. Furthermore, the correspondence 
$w\mapsto \psi_w$ gives us a map $\psi: W \to \nn$. 

\m It is easy to see that the bordism class of $\bigl [ W,\psi]$ 
depends only on the bordism classes 
$\ga_1, \ga_2$. 

\begin{thm}\label{bracket} 
The operation $[-,-]: \oma_i (\NN)\otimes \oma _j(\NN) \to \oma 
_{i+j+2n}(\NN )$, $i,j\in \Z,$ 
converts the abelian group $\oma _* (\NN )$ into the graded Lie algebra, 
i.e. the 
operation $[-, -]$ has the following properties: 

\begin{description} 
\item[1] $[\ga_1, \ga_2]=(-1)^{(|\ga_1|+n)(|\ga_2|+n)} [\ga_2, 
\ga_1]$, for all $\ga_1, 
\ga_2 \in \oma_ *(\NN)$; 
\item[2] $(-1)^?[[\ga_1, \ga_2], \ga _3]+(-1)^?[[\ga_2, \ga_3], 
\ga_1]+(-1)^?[[\ga_3, \ga_1], \ga_2]=0$, for all $\ga_1, \ga_2, 
\ga_3\in \oma_*(\NN)$; 
\item[3] $[\ga_1, \ga_2 + \ga _3]=[\ga_1, \ga _2] +[\ga_1, 
\ga_3]$, for all $\ga_1, \ga_2, \ga_3\in \oma_*(\NN)$. 
\end{description} 
\end{thm}

\pp We are not able to put proper signs in the proposed Jakobi identity 
(2). (Both authors are sure that such identity exists, but up to present 
the authors get different answers.) However, we can prove the Jacobi 
identity modulo 2, and we do it now.
Without loss of generality we assume that that the bordisms $\ga_1, 
\ga_2, \ga_3$ have adjoint maps of the form 
$\ov\ga_i:F_i\times P_i\to M$, $i=1,2,3$, with $F_i$ connected.
Then all the garlands that appear when one opens the Lie brackets in the 
Jacobi identity consist of  $P_1$, $P_2$, and $P_3$; and all garlands 
have two two-point marks 
that attach two out of three manifolds $P_1$, $P_2$ and $P_3$ to the 
third one. The appearing types of garlands subdivide in three different 
types depending on which one of $P_1$, $P_2$ and $P_3$ contains two out 
of four points of two two-point marks. (This is the manifold to which 
the other two are attached.)

The garlands where two out of four double points are on $P_1$ appear 
from the first and the third brackets of the Jacobi identity. 
The manifolds $S_1$ and $S_3$ describing the corresponding bordisms in 
$\NN$ coming from the first and the third bracket are described by the 
following set-theoretic conditions.
$$
S_1=\bigl \{ (f_1, n_1^1, f_2, n_2, n_1^2, f_3, n_3)\bigr \}\subset 
F_1\times N_1\times F_2\times N_2\times N_1\times F_3\times N_3
$$ 
such that $\ov \ga_1 (f_1, n_1^1)=\ov \ga_2 (f_2, n_2)$  and  $\ov 
\ga_1(f_1, n_1^2)=\ov \ga_3 (f_3, n_3).$

$$
S_3=\bigl \{ (f_3, n_3, f_1, n_1^2, n_1, f_2, n_2)\bigr \}\subset 
F_3\times N_3\times F_1\times N_1\times N_1\times F_2\times N_2
$$
such that  $\ov \ga_3(f_3, n_3)=\ov \ga_1(f_1, n_1^2)$ and 
$\ov \ga_1(f_1, n_1^1)=\ov \ga_2 (f_2, n_2).$

\m
It is clear that $S_1$ and $S_2$ can be transformed to each other via a 
coordinate permutation. Since the mappings of the corresponding garlands 
are determined by the values of $f_1, f_2, f_3$, we get that the 
mappings of garlands appearing from the corresponding points of $S_1$ 
and $S_3$ are the same. Thus this type of garlands in the Jacobi 
identity cancels out when one considers bordisms with $\Z_2$ 
coefficients.

\m
The proof of the fact that the other two types of garlands appearing 
from the Jacobi identity cancel out is the same. The proof in the case 
where the garlands corresponding to $\ga_1$, $\ga_2$, and $\ga_3$ 
consist of more than one copy of the manifold $P$ is obtained in the 
same way since the Lie bracket $[-,-]$ is defined via pull-backs over 
pairs of participating copies of $P$.\qed

\section{The homomorphisms $\li$, $\pr$, and $\Delta$}\label{lip} 

Consider a map $\gf: F\to\NN$ with $F$ connected, and let $\ov\gf: F 
\times N\to M$ be the adjoint marked map. Given a point $(f,n)\in 
F\times N$, we define a map $\psi_{f,n}: (N,\fm) \to M$ as follows: 
$\fm=\{n\}$ and has grading 1, and $\psi_{f,n}(n')=\ov\gf(f,n')$ for 
all $n'\in N$. All the other marks increase their grading by 1. So, 
$\psi_{f,n}\in \NN$, and we define $\psi: F\times N \to \NN$ as 
$\psi(f,n)=\psi_{f,n}$. If we have a map $\gf: F \to 
M=\nn(0;0)\subset \nn$, we define $\psi: F \times P \to M$ to be the 
composition of $\gf$ and the projection $F \times P \to P$. 

Now, the correspondence $\gf\mapsto \psi$ yields a homomorphism 
\begin{equation}\label{lift} 
\li: \oma_i(\NN) \to \oma_{i+n}(\NN). 
\end{equation} 

\begin{defin} we define the homomorphism $\pr: \oma_i(\nn) \to 
\oma_i(\nn)$ as follows. Given a map $\ov\ga: F \to \NN$ with the 
adjoint map $\omega: F \times N \to N$, the map $\pr$ just changes the 
marks on $N$ as follows: 
\begin{description} 
\item[1] if $g(\fm)\ne 1$ then $\pr$ decreases grading by 1; 
\item[2] if $g(\fm)=1$ and $|\fm|>1$ then $\pr$ increases grading by 
1; 

\item[3] if $g(\fm)=1$ and $|\fm|=1$ then $\pr$ erases the mark. 
\end{description} 
\end{defin} 

Using the definitions of $\pr, \li, \bullet$ and $[-,-]$ one gets the 
following result.

\begin{prop} For all $\ga_1,\ga_2\in \bor_*(\NN)$ we have 
$$ 
[\ga_1,\ga_2]=\pr(\li(\ga_1)\bullet\li(\ga_2)). 
$$ 
\end{prop} 

\begin{prop}\label{prli=0} 
If the manifold $P$ is a boundary, then $\pr\circ\li\equiv 0$. 
\end{prop} 

\pp Let $\ga\in \oma_i(\nn)$ be represented by a map $\ov\ga: F \to 
\NN$. It is easy to see that the element $\pr\circ\li(\ga)$ is 
represented by the map 
$$ 
\CD 
F\times N @>\text{projection}>> F @>\omega >> \NN. 
\endCD 
$$ 
But $N$ is boundary since $P$ is, and thus the above map is 
zero-bordant. 
\qed 

\begin{defin} We define the homomorphism $\Delta : \oma_i(\NN)\to 
\oma_{i+n}(\NN)$ by setting $\Delta=\li\circ \pr$.
\end{defin}

\begin{cor}\label{Delta}
If the manifold $P$ is a boundary, then $\Delta^2=0$.
\end{cor}

\pp This follows from \propref{prli=0}, because 
$$
(\li\circ \pr)\circ (\li \circ \pr)=\li \circ (\pr \circ \li ) \circ 
\pr=0.
$$
\qed

\section{Batalin--Vilkovyski and Gerstenhaber structures} 

Given a graded commutative algebra $A=(A,+,\bullet)$, consider an 
additive homogeneous homomorphism $\Delta: A_*\to A_{*+n}, n>0$. A 
quadruple $A, +,\bullet, \Delta)$ is called a {\it Batalin--Vilkovyski 
algebra of degree $n$} if $\Delta^2=0$ and (cf. Getzler 
\cite{Getzler}). 
\begin{equation}\label{getzler} 
\begin{aligned} 
\Delta(a\bullet b \bullet c) &=\Delta(a\bullet b)\bullet c 
+(-1)^{|a|n}a\bullet \Delta(b\bullet c)\\ 
&+(-1)^{(|a|+n)|b|}b\bullet \Delta(a\bullet c) -\Delta(a)\bullet b 
\bullet c\\ &-(-1)^{|a|n}a \bullet \Delta (b)\bullet c 
-(-1)^{n(|a|+|b|)}a\bullet b \bullet \Delta(c) 
\end{aligned} 
\end{equation} 

\m A {\it graded Gerstenhaber algebra of degree $n$} (called also {\it 
Poisson algebra} or {\it braid algebra} of degree $n$) is a graded 
commutative algebra $V$together with an operation $\{-,-\}: V \otimes V 
\to V$ satisfying the following relations (cf.  
Gerstenhaber~\cite{Gerstenhaber}): 
\begin{equation}\label{poissoncommut} 
\{a,b\}=(-1)^n (-1)^{(|a|+n)(|b|+n)}\{b,a\} 
\end{equation} 
\begin{equation}\label{poissonjacobi} 
\{a,\{b,c\}\}=\{\{a,b\},c\}+(-1)^{(|a|+n)(|b|+n)}\{b,\{a,c\}\}, 
\end{equation} 
and
\begin{equation}\label{gerst} 
\{a, b\bullet c\}=\{a, b\}\bullet 
c+(-1)^{(|a|+n)b}b\bullet\{a, c\} 
\end{equation} 
\begin{prop} Let $(A, \bullet,\Delta)$ be a graded Batalin--Vilkovyski 
algebra. Set 
$$ 
\{a,b\}=(-1)^{|a|n}\Delta(a\bullet b) - (-1)^{|a|n}\Delta(a) \bullet b 
-a\bullet \Delta(b) 
$$ 

Then the quadruple $(A, +,\bullet,\{-,-\})$ is a graded Gerstenhaber 
algebra of degree $n$. 

Furthermore, the Leibnitz rule 
\begin{equation}\label{leibnitz} 
\Delta\{a, b\}=(-1)^{n+1}\left(\{\Delta a, b\}+(-1)^{|a|n+1}\{a, 
\Delta b\}\right) 
\end{equation} 
holds. 
\end{prop} 

\pp Getzler~\cite[Prop. 1.2] {Getzler} proved these equalities for 
$n=1$ (cf. also Penkava--Schwarz~\cite{PenkavaSchwarz}). His approach 
works for arbitrary $n$, since we properly defined the signs as 
functions of $n$. 
\qed

We hope but cannot prove in general that $\Delta$ satisfies relation 
\eqref{getzler}. However, we have several examples where $\Delta$ 
defined as in Corollary~\ref{Delta}, does satisfy~\eqref{getzler}, and 
thus $\oma_*(\nn)$ possesses the structures of Batalin--Vilkovyski and 
Gerstenhaber algebras.

\section{Comparison with Chas--Sullivan string theory}

Here we just notice that for $P=S^1$ our theory is parallel to 
Chas-Sullivan string theory~\cite{ChasSullivan}. 

\m In \cite{ChasSullivan} Chas and Sullivan consider the space 
$LM=\Lambda M/SO(2)$ where $\Lambda M$ is the space of smooth free 
loops on $M$ (i.e., maps $S^1 \to M$) and the $SO(2)$-action on $\Lambda 
M$ is induced by the tautological $SO(2)$-action on $S^1=SO(2)$. 
Furthermore, they define operations $\bullet$ and $[-,-]$ on $H_*(LM)$. 
Certainly, these operations can also be defined for bordisms instead of 
homology, and we get the operations 
$$ 
\bullet :\oma_i(LM)\otimes \oma_j(L M) \to \oma_{i+j}(LM) 
$$and 
$$ 
[-,-]: \oma_i(LM)\otimes \oma_j(LM) \to \oma_{i+j+2}(LM). 
$$ 

\begin{rem} In fact, we prefer to work with bordisms because it is a 
more geometric theory, it is easy to control transversality, etc. It 
seems, however, that our approach works for homology also, if we will be 
able to work properly with transversality. 
\end{rem}

\m Some troubles appear if one tries to generalize Chas--Sullivan 
construction for an arbitrary manifold $P$ instead of $S^1$. There are 
several of them, but we discuss the following which looks most 
important. Chas--Sullivan use the co-$H$-space structure on $S^1$; this 
allow them to regard, say, a wedge (composition) of loops as a loop. 
However, in the class of closed manifold only homotopy spheres are 
co-$H$-spaces. So, if we want to consider any $P$ instead of $S^1$ then 
we should, roughly speaking, convert $P$ into a co-$H$-space. 

\m It is instructive to consider an algebraic analog of the situation. 
If we have an abelian group $A$ and want ``to convert it to a ring'', we 
consider the tensor algebra $T(A)$. So, conceptually, for a manifold $P$ 
a role of an analog of the tensor algebra can play a space $P\sqcup 
P\vee P \sqcup \ldots$ (the union of garlands). This is not a manifold, 
but we can resolve the singular points of wedges via considering the 
marked maps. This ideas lead us to spaces $\NN$, which can be considered 
as analogs of $\Lambda M$. It is also worthy to mention that 
Chas--Sullivan uses $SO(2)$-action in order to work with wedges of 
non-pointed spaces. For arbitrary $P$, this problem is solved by using 
marked maps. 

\m We finish this notes with a vague remark that a map of connected 
garland based on $S^1$ can be regarded as a map of circle. We do not 
dwell and explain here in detail the parallelism between our theory (for 
$P=S^1$) and Chas--Sullivan one,  but we hope to do it somewhere else 
(or in later version of this manuscript).


\begin{thebibliography}{99999} 

\bibitem{Arnoldcurves} 
V.~I.~Arnold: {\it Plane curves, their invariants, perestroikas 
and 
classifications, Singularities and Bifurcations} 
(V.I.~Arnold, ed.) 
Adv. Sov. Math., Vol. 21 (1994), pp. 39--91 

\bibitem{Arnold} 
V.~I.~Arnold: {Invariants and perestroikas of wave fronts on 
the plane, Singularities of smooth maps with additional 
structures}, 
Proc. Steklov Inst. Math., Vol. 209 (1995), pp. 11-- 
64. 

\bibitem{Arnoldbook} 
V.~I.~Arnold: {\em Mathematical Methods of Classical 
Mechanics,} Second 
edition, Graduate Texts in Mathematics, {\bf 60}, Springer- 
Verlag, 
New-York (1989) 

\bibitem{ChasSullivan} 
M.~Chas and D.~Sullivan: {\em String Topology,\/}, 
math.GT/9911159 at the http://xxx.lanl.gov 

\bibitem{ChernovRudyak} 
V.~Chernov (Tchernov) and Yu.~B.~Rudyak: {\em Affine 
linking numbers and 
causality relations for wave fronts,\/} submitted for 
publication (2002) also available as a preprint 
math.GT/0207219 at the http://xxx.lanl.gov 

\bibitem{ChernovRudyakWinding} 
V.~Chernov (Tchernov) and Yu.~B.~Rudyak: {\em Affine 
winding numbers and front propagation,\/} submitted for 
publication (2002) also available as a preprint 
math.GT/0301117 at the http://xxx.lanl.gov 

\bibitem{Gerstenhaber} 
M. Gerstenhaber,{\em The cohomology structure of an associative ring,\/} 
Ann. of Math. (2) 78 (1963) 267--288. 


\bibitem{Getzler} 
E. Getzler: {\em Batalin--Vilkovyski Algebras and Two-Dimensional 
Topological Field Theories} Comm. Math. Phys. {\bf 159} (1994), pp. 
265--285. 

\bibitem{Goussarov1} 
M.~Goussarov: A talk on 12.1.1987 where the finite order invariants 
were first defined 
at the Rokhlin seminar at the Steklov Institute for Mathematics, St 
Petersburg, 
Russia. 


\bibitem{Goussarov2} 
M.~Gusarov: 
{\em A new form of the Conway-Jones polynomial of oriented links.\/} 
(Russian. 
English summary) 
Zap. Nauchn. Sem. Leningrad. Otdel. Mat. Inst. Steklov. (LOMI) 193 
(1991), Geom. 
i Topol. 1, 4--9, 161. 

\bibitem{PenkavaSchwarz} 
M.~Penkava, A.~Schwarz,{\it: On some algebraic structures arising in 
string theory}. Perspectives in mathematical physics, 
219--227, Conf. Proc. Lecture Notes Math. Phys., III, Internat. Press, 
Cambridge, MA, 1994. 

\bibitem{Rudyak}
Yu. B. Rudyak: {\em On Thom Spectra, Orientability, and Cobordism}, 
Springer,Berlin Heidelberg New York, 1998

\bibitem{Spanier} 
E.~Spanier: {\em Algebraic topology.\/} 
Corrected reprint of the 1966 original. Springer-Verlag, New York, 1996 

\bibitem{Stong} 
R.~Stong: {\it Notes on cobordism theory}, Mathematical notes, Princeton 

University Press, Princeton, N.J. 1968 

\bibitem{Switzer} 
R.~Switzer: {\it Algebraic topology---homotopy and homology}, Die 
Grundlehren der mathematischen Wissenschaften, Band 212. Springer, 
Berlin Heidelberg New York, 1975. 

\bibitem{Thom} 
R.~Thom: {\em Quelques propri\'et\'es globales des vari\'et\'es 
diff\'erentiables,\/} 
Comment. Math. Helv. 28, (1954). pp. 17--86. 

\bibitem{Vassiliev} 
V.~A.~Vassiliev: {\em Cohomology of knot spaces,\/} 
Theory of singularities and its applications, pp. 23--69, 
Adv. Soviet Math., 1, 
Amer. Math. Soc., Providence, RI, (1990) 

\end{thebibliography}
\end{document}